\documentclass{amsart}
 \usepackage{amsmath}
\usepackage{amscd}
\usepackage{latexsym}
\usepackage{amssymb}
\usepackage{mathtools}
\usepackage{xcolor}
\newtheorem{cor}{Corollary}

\begin{document}
\title{On the quantization of $C^{\infty}(\mathbb R^d)$}
\author{Murray Gerstenhaber}
\thanks{The author thanks Jim Stasheff for his help in the preparation of this paper.}
\address{Department of Mathematics\\  University of Pennsylvania\\
  Philadelphia, PA 19104-6395, USA}
\email{mgersten@math.upenn.edu}

\begin{abstract}\noindent 
The Basic Universal Deformation Formula is proven and applied to show  that Weyl algebras, which encode Heisenberg's uncertainty principle, are effective deformations of polynomial rings, and that uncertainty is necessary for stability. Deformation problems may have associated modular groups an algebraic example of which is given. Poisson structures on $C^{\infty}(\mathbf R^d)$, 
shown by Kontsevich to be infinitesimal deformations integrable to full deformations, here are  
 shown to be the skew forms of those infinitesimal deformations of $C^{\infty}(\mathbf R^d)$ with vanishing primary obstructions.  
In dimension 3 any smooth multiple of such an infinitesimal again has vanishing primary obstruction. This exceptional property suggests that in our universe a large  disturbance like the Big Bang can be confined to an arbitrarily small interval in time and almost completely confined to an arbitrarily small region in space.
\end{abstract}
 
\dedicatory{To the memory of my beloved wife, Ruth P. Zager M.D., z''l}
\maketitle
\newtheorem{theorem}{Theorem}
\newtheorem{corollary}{Corollary}
\newtheorem{lemma}{Lemma}
\newtheorem{defn}{Definition}

\newcommand{\R}{\ensuremath{\mathbb{R}}}
\newcommand{\C}{\ensuremath{\mathbb{C}}}
\newcommand{\Cinf}{\ensuremath{C^{\infty}(\mathbb{R}^d)}}
\newcommand{\Q}{\ensuremath{\mathbb{Q}}}
\newcommand{\Z}{\ensuremath{\mathbb{Z}}}
\newcommand{\cP}{\ensuremath{\mathcal{P}}}
\newcommand{\cM}{\ensuremath{\mathcal{M}}}       
\newcommand{\cI}{\ensuremath{\mathcal{I}}}
\newcommand{\ci}{\ensuremath{\mathcal{i}}}
\newcommand{\cK}{\ensuremath{\mathcal{K}}}
\newcommand{\cL}{\ensuremath{\mathcal{L}}}
\newcommand{\cW}{\ensuremath{\mathcal{W}}}
\newcommand{\bM}{\ensuremath{\mathbb{M}}}
\newcommand{\bS}{\ensuremath{\mathbb{S}}}
\newcommand{\cX}{\ensuremath{\mathcal{X}}}
\newcommand{\cA}{\ensuremath{\mathcal{A}}}
\newcommand{\cD}{\ensuremath{\mathcal{D}}}
\newcommand{\cR}{\ensuremath{\mathcal{R}}}
\newcommand{\bA}{\ensuremath{\mathbb{A}}}
\newcommand{\bR}{\ensuremath{\mathbb{R}}}
\newcommand{\cB}{\ensuremath{\mathcal{B}}}
\newcommand{\cS}{\ensuremath{\mathcal{S}}}
\newcommand{\pr}{\ensuremath{\preceq}}
\newcommand{\op}{\ensuremath{\mathrm{op}}}
\newcommand{\bfk}{\ensuremath{\mathbf{k}}}
\newcommand{\ctD}{\ensuremath{\overline{\mathcal{D}}}}

\newcommand{\px}{\ensuremath{\partial_x}}
\newcommand{\py}{\ensuremath{\partial_y}}
\newcommand{\pz}{\ensuremath{\partial_z}}
\newcommand{\pw}{\ensuremath{\partial_w}}
\newcommand{\pl}{\ensuremath{\partial_}}

\newcommand{\btr}{\ensuremath{\blacktriangleright}}
\newcommand{\vtr}{\ensuremath{\vartriangleright}}
\newcommand{\tR}{\ensuremath{\overline{R}}}
\newcommand{\car}{\ensuremath{\circlearrowright}}
\newcommand{\g}{\gamma}
\newcommand{\la}{\lambda}
\newcommand{\om}{\omega}
\newcommand{\Om}{\Omega\newcommand{\tR}{\ensuremath{\overline{R}}}}
\newcommand{\G}{\Gamma}
\newcommand{\z}{\zeta}
\newcommand{\Rp}{{R}^+}
\newcommand{\Aut}{\operatorname{Aut}}
\newcommand{\Der}{\operatorname{Der}}
\newcommand{\iso}{\operatorname{iso}}
\newcommand{\dr}{\operatorname{dR}}
\newcommand{\sk}{\operatorname{sk}}
\newcommand{\id}{\operatorname{id}}
\newcommand{\pex}{\operatorname{pex}}
\newcommand{\frg}{{g}}
\newcommand{\lbrr}{\linebreak\raggedright}

\newcommand{\cpr}{\ensuremath{\!\!\smile\!\!}}

\newcommand{\U}{\Upsilon}
\vspace{-7mm}
\date{}

\section{History} \label{History}
Algebra deformation theory has a long background in the deformation theory of complex analytic structures but in its present form originated in the author's papers \cite{G:Cohomology}, \cite{G:Def} and was further developed in \cite{G:Def2}, \cite{G:Def3}, \cite{G:Def4}. Riemann studied the `moduli' on which compact Riemann surfaces depend which may be considered  as deforming when the moduli are varied. When different points in the parameter space yield structures isomorphic to the original there is a modular group operating on that space. An algebraic example is given here, analogous to the classical modular group associated with deformations of compact Riemann surfaces of genus one. Infinitesimal deformations of Riemann surfaces were first formalized by Teichm{\"u}ller, \cite{Teich:Extremale}  as quadratic differential but the break-through definition of infinitesimal deformations for a complex manifold of arbitrary dimension 
 as elements of the cohomology in dimension one of the manifold with coefficients  in its sheaf of germs of holomorphic tangent vectors 
 is due to Fr{\"o}licher and Nijenhuis, \cite{FN}. This led to the comprehensive  works of  Kodaira and Spencer; for an overview see \cite{Kodaira:Complex}.  However, obstructions as well as  jump deformations, \cite{Coffee:Filtered}, \cite{G:Def4},\cite[\S 7, p. 20ff]{GS:AlgebraicCohomology}, were not clearly understood until the development of algebraic deformation theory. These are discussed in \S\ref{jump}.
Deformation Quantization, an extraordinary advance in algebraic deformation theory introduced by Bayen, Flato, Fr{\o}nsdal, Lichnerowitz, and Sternheimer, \cite{BFFLS}, allows quantum theory to be understood and the spectrum of the hydrogen atom to be computed without the use of wave mechanics or Schr{\"o}dinger's equation. For a summary and some later developments, see \cite{Sternheimer:20YearsAfter} and \cite{DitoSternheimer:DefQuant}. Subsequent to \cite{BFFLS}, in the context of algebras,``quantization" has become become almost synonymous with "deformation" and is so used here.

\section{Overview}\label{Overview}
 Algebraic deformation theory asks the following question: Given an algebra $\cA$ over some commutative unital ring $\bfk$, in what ways can one create a ``deformation" of that  algebra with  multiplication of the form
\begin{equation}\label{star}
a\star b = ab + \hbar F_1(a,b) + \hbar^2F_2(a,b ) + \cdots \quad
\end{equation} 
while remaining in the same ``equationally defined category", so called in \cite{G:Def2}, or in present terminology, category of algebras defined over a particular operad. The notation and term ``star product''  were introduced in \cite{BFFLS}. The symbol $\hbar$ originally denoted the reduced Planck's constant, $h/2\pi$, 
where $h \approx 6.626176 \times 10^{-34}$ joule-seconds is Planck's original constant. Its use to denote a deformation parameter is now common and derives from the fact that Planck's constant may, in fact, be viewed as such, cf \cite{BFFLS}, \cite{G:PathAlgebras}.\footnote{Correction to Theorem 8 of \cite{G:PathAlgebras}: In view of the remark at the end of \S2 (\emph{not} \S1) it should read, ``...the coherent twist of the path algebra induced by $\omega$ is trivial \emph{only if} the class of $\omega$  is trivial as an element of $H^2(\mathcal M, \R/\tau\R)$."}

For associative algebras, the $F_i$ in \eqref{star} are Hochschild 2-cochains of $\cA$ with coefficients in $\cA$ itself, tacitly extended to be defined over $\cA[[\hbar]]$; we may write $F_0$ for the original multiplication in $\cA$.
 Here $F_1$ must be a 2-cocycle, usually called the infinitesimal of the deformation and a main problem is to determine which infinitesimals give rise to full deformations.  Note that if $D_1, D_2$ are derivations then $D_1 \smile D_2$ is a 2-cocycle.
 Given a 2-cocycle $F_1$, there generally is a sequence of cohomological obstructions to constructing the necessary $F_i, i \ge 2$. The analogue of the foregoing in the Kodaira-Spencer approach to the deformation of complex analytic manifolds is already present in an insufficiently recognized paper by Douady \cite{Douady:Obstruction}. A full deformation as in \eqref{star}, if it exists, is often called an integral of $F_1$, something which depends only on the cohomology class of $F_1$. While it is really the class which is  the infinitesimal, of which  $F_1$ is only a representative, it will be convenient also to call $F_1$ `the' infinitesimal of the deformation. 
 
 A {\em Universal Deformation Formula} (UDF) exhibits in closed form an explicit integral for some class of 2-cocycles, considered as infinitesimal deformations. The first such, called here the {\em Basic} UDF, was introduced  in \cite[Ch. II, Lemma 1, p. 13]{G:Def3} and stated there for ``composition complexes'', amongst which are the cohomology groups of algebras, coalgebras, and simplicial complexes, but the first proof appears in \cite{GiaqZhang:UDF}. While needed here only for deformations of  associative algebras the proof in \S\ref{BasicUDF} holds for all composition complexes.
 
 The Basic UDF asserts that if $D_1, D_2$ are commuting derivations of an algebra $\cA$ then a suitably defined exponential of $D_1\cpr D_2$ is an integral. 
 The star product of \eqref{star} is, however, only a formal power series, while what is generally wanted is that the deformed algebra be defined over the original coefficient ring $\bfk$. It may not be meaningful, \textit{a priori}, to specialize the deformation parameter $\hbar$ to an element of $\bfk$.  When this is so for a particular UDF operating on some algebra, the operation will be called effective.  In \S\ref{Weyl} it is shown, using the Basic UDF, that the   
 first Weyl algebra, which encapsulates Heisenberg's uncertainty relation between position and momentum, is an effective deformation of the polynomial ring in two variables; although \eqref{star} is formally a power series, when applied to any pair of terms in the original algebra it terminates in a finite number of steps. Clifford algebras are, likewise, effective deformations of graded polynomial rings, when degree is taken into account. The Basic UDF, \S\ref{non-com}, also gives another approach to the first UDF with non-commuting derivations, \cite{CGG:noncomm}.
 
 Kontsevich's UDF, \cite{Kontsevich:Poisson}, provides an explicit, but again purely formal, integral for any Poisson infinitesimal deformation  in the form of a power series in partial derivatives. 
It uses the fundamental works of Stasheff, \cite{Stasheff:Hspaces}, and Schlessinger and Stasheff, \cite{SchlessStash}, which introduced the concepts of homotopy associative algebras and homotopy Lie algebras, respectively. 
The concept of algebras up to homotopy is now also understood for algebras over operads. 
 Tamarkin, \cite{Tamarkin:formality}, almost simultaneously gave another proof of the integrability of Poisson structures, but by quite different methods. Cattaneo, Felder, and Tomassini, \cite{CattFeldTom}, extended these results to the quantization of Poisson structures on a manifold. 
The coefficients in Kontevich's UDF are not rational and are difficult to compute. Dolgushev, \cite{Dolg:rational}, later showed that if a solution exists, then there is a rational one. Determining the Poisson structures on $\Cinf$ remains an open problem.
A power series integral to an infinitesimal deformation generally does not converge when the infinitesimal is smooth but not analytic. However, effective quantization is sometimes possible by using Fourier series,  \S \ref{Fourier}.

An infinitesimal deformation, viewed as a cohomology class, has a unique representative 2-cocycle which is antisymmetric or skew. The primary obstruction to integrating an infinitesimal deformation of a commutative algebra $\cA$ in which 6 is a unit is shown to vanish if and only if its unique  skew representative is Poisson.  Any infinitesimal deformation of $\Cinf$ with vanishing primary obstruction is therefore, by the work of Kontsevich, integrable to a formal deformation.  It is shown also that the primary obstruction to an infinitesimal deformation of $\Cinf$ vanishes if and only if the coefficients of a representative cocycle satisfy a system of $\binom{d}{3}$ non-linear partial differential obstruction equations. 
 When $d = 2$ there are no obstructions; all infinitesimal deformations give rise to quantizations. For $d=3$, the single obstruction equation has a property not present in higher dimensions; it allows an arbitrarily large deformation to be localized to an arbitrarily small region in space and an arbitrarily small interval of time, \S10. This may help us to understand the Big Bang.

 \section{The Basic Universal Deformation Formula}  \label{BasicUDF}								
 Throughout, $\bfk$ will denote a commutative unital ground ring and $\cA$ will 
 be a $\bfk$-algebra. When we have a deformation of $\cA$ with star product as in \eqref{star}, note that if $a,b \in\cA$ then $a\star b$ lies not in $\cA$ but in $\cA[[\hbar]]$, the ring of power series in $\hbar$ with coefficients in $\cA$. While in general it is not possible to specialize $\hbar$ to an element of the ground ring $\bfk$, it will be shown that this is possible in some important cases.

The condition that star products must remain associative can be expressed
 using special cases of the composition products $\circ_i$ and $\circ$ introduced in \cite{G:Cohomology}. If $F,G$ are 2-cochains of $\cA$ then $F\circ_1G,\,F\circ_2G, \, F\circ G$ are the 3-cochains defined by setting, respectively,  
\begin{equation}\label{circ12}
(F\circ_1G)(a,b,c) = F(G(a,b),c), \quad(F\circ_2G)(a,b,c) = F(a,G(b,c)),
\end{equation}
and 
\begin{equation}\label{circ}
(F\circ G)(a,b,c) = (F\circ_1G-F\circ_2G)(a,b,c) = F(G(a,b),c) - F(a, G(b,c)).
\end{equation} 
 With these, the requirement that $a\star (b\star c) - (a\star b) \star c$ = 0 can be expressed by the equations
\begin{equation}\label{basic}
\sum_{i= 0}^n
F_i\circ F_{n-i}\quad = 0
\quad \text{for all} \quad n.
\end{equation}
Transposing to the right those terms where either $i=0$ or $i= n$, \eqref{basic} can be rewritten as
\begin{equation}\label{Basic2}
\sum_{i+j=n, \, i,j > 0}
F_i\circ F_j = -\delta F_n,
\end{equation}
where $\delta$ is the Hochschild coboundary operator. 
 When $n = 1$ the left side vanishes, showing that  $F_1$ is, as mentioned, a 2-cocycle.  
 
 Star products $a\star  b = \sum \hbar^iF_i$ and $a\star'  b = \sum \hbar^iF_i'$ are called  
 gauge equivalent if there is 
a one-parameter family $\gamma$ of $\bfk$-linear automorphisms of $\cA$ of the form 
$\gamma(a) = a+ \hbar \gamma_1(a) +  \hbar^2 \gamma_2(a_) + \cdots$, where the $\gamma_i$ are linear maps from $\cA$ to itself, again tacitly extended to be defined over $\cA[[\hbar]]$, such that $a\star'b = \gamma^{-1}(\gamma(a)\star\gamma(b))$. The deformed algebras that $\star$ and $\star'$ define on the underlying $\bfk$-space of $\cA$ are then isomorphic, and one has $F_1' = F_1 +\delta\gamma_1$.  As $\gamma_1$ can be arbitrary, $F_1$ can be replaced in a gauge equivalent deformation by any cohomologous 2-cocycle, showing as remarked, that the integrability of $F_1$ depends only on its cohomology class.  

When $F_1$ is a cocycle  then so is $F_1\circ F_1$. Its cohomology class in $H^3(\cA,\cA)$, which depends only on the class of $F_1$, is its primary obstruction, a term also commonly applied to $F_1\circ F_1$ itself. If $F_1\circ F_1$ is a coboundary then one can choose an $F_2$ with $-\delta F_2 = F_1\circ F_1$ and one can ask if an $F_3$  exists so that one can continue building the series \eqref{star}. However, \eqref{Basic2} with $n =3$ shows that one may encounter another obstruction in $H^3(\cA,\cA)$, and so on indefinitely. If, with a given 2-cocycle $F_1$, we are able to  construct a series such as that in \eqref{star}, then that integral is also said to quantize $\cA$. When $\cA$ is the algebra $\Cinf$ of smooth functions on $\R^d$ the integral is also said to quantize $\R^d$.
 While in principle one may encounter an infinite sequence of obstructions,  
 Kontsevich's work,  \cite{Kontsevich:Poisson}, together with what will be shown here, implies that an infinitesimal deformation $F_1$ of $\Cinf$ has only the primary obstruction; if that one obstruction vanishes, then there is a star product with the given $F_1$.  
 
 An algebra $\cA$ for which all deformations are gauge equivalent to the trivial deformation, i.e., the one where the star product is just the original multiplication, is called rigid. This will be the case if every infinitesimal deformation is ultimately obstructed and certainly if $H^2(A) = 0$, in which case the algebra is called absolutely rigid or, preferably, stable in the terminology of \cite{FN}. A tensor product of stable algebras need not be stable, as will be seen.

When $H^3(\cA,\cA) = 0$, every infinitesimal deformation is integrable.  In particular, this is the case for $C^{\infty}(\R^2)$ by the Hochschild-Kostant-Rosenberg (HKR) Theorem, \cite{HKR}.  That theorem asserts that $H^*(\Cinf)$ is  isomorphic, as a module over $\Cinf$, to the exterior algebra generated by the partial derivatives, $\pl1, \dots ,\pl d$ with respect to the coordinates $x_1,\dots, x_d$.  In particular, when $d=2$, one has $H^3(C^{\infty}(\R^2)) = 0$, so there are no obstructions to any infinitesimal deformation of  $C^{\infty}(\R^2)$.  (The original formulation of the HKR Theorem was more restrictive, for a proof of the present form, see  Roger, \cite{Roger}.) Note, however, that the module isomorphism does not carry the exterior product on $\wedge^* \Der{\cA}$ to the cup product on $H^*(\cA,\cA)$, cf \cite{CalaqueVandenBergh}.

If $D_1,\dots, D_n$ is a sequence of derivations of $\cA$ (which need not commute or be distinct), and $I = (i_1,\dots,i_r)$ is a subsequence of $(1,\dots, d)$, set $D_I = D_{i_1}D_{i_2}\cdots D_{i_r}$ and let $I^c$ denote the ordered complement of $I$. Then $\delta(D_1D_2\cdots D_n) = -\sum_ID_I\cpr D_{I^c}$, where the sum is over all ordered proper subsequences of $(i_1,\dots, i_r)$, i.e., neither empty nor the whole (a sum over ``unshuffles").
In particular, for any pair of derivations $D_1,D_2$ of $\cA$, we have $\delta(D_1D_2) = - (D_1\cpr D_2 + D_2\cpr D_1)$. Therefore, $D_1\cpr D_2$ is cohomologous to $-D_2\cpr D_1$ and also to $(1/2)(D_1\cpr D_2 - D_2\cpr D_1)$, provided that 2 is a unit.   It follows from the HKR Theorem
 that every 2-cocycle of $\Cinf$ is cohomologous both to one in normal form $\sum_{1\le i < j \le d}a_{ij}\pl i\cpr\pl j$ and to one in skew form $\sum_{1\le i < j \le d}(a_{ij}/2)(\pl i\cpr\pl j - \pl j\cpr\pl i)$. 

All Hochschild cohomology groups of an algebra $\mathcal A$ considered here will, unless stated otherwise, have coefficients in $\mathcal A$ itself as a bimodule. That cohomology is frequently called the regular (Hochschild) cohomology and will be denoted simply $H^*(\mathcal A)$, and similarly for cochains. Note in what follows that if $a$ is a central element of $\cA$ and $F$ a cocycle of any dimension, then $aF$ is again a cocycle and if $a$ is central and  $D$ is a derivation, then $Da$ is again central. 

The product of $1$-cochains of an associative algebra $\cA$ is always well defined as their composition. For $2$-cochains of the form $f\cpr g$, where $f,g$ are $1$-cochains, one can not generally define a product by setting $(f_1\cpr g_1)(f_2 \cpr g_2) = f_1f_2 \cpr g_1g_2$, for if $a$ is a central element of $\cA$, then as $2$-cochains one has $a(f\cpr g) = af \cpr g = f\cpr ag$ but such changes in representation will usually change the product. Suppose, however, that we have a set of commuting derivations $\{D_i\}$ of $\cA$.  
When $I = \{i_1, \dots, i_r\}$ is a set of indices (repetitions allowed) of the $D_i$, set, as before, $D_I = D_{i_1}\cdots D_{i_r}$, and similarly for some $D_J$. As the derivations now are assumed to commute, the ordering does not matter. 
(When the set $I$ is empty, interpret $D_I$ as the identity map $\cA$ and similarly for $D_J$; when $I$ and $J$ are both empty then $D_I \cpr D_J$ will denote the multiplication map.) The foregoing problem then does not arise when multiplication is restricted to $2$-cochains of the form $D_I \cpr D_J$.  In particular, if $D_1, D_2$ are commuting derivations of $\cA$ then $(D_1 \cpr D_2)^n = D_1^n \cpr D_2^n$ is well defined, as is $\exp (D_1\cpr D_2) = \sum_{n=0}^{\infty}\frac{1}{n!}(D_1 \cpr D_2)^n$  when $\cA$ is defined over $\Q$.

In the following, `formal' indicates that no assertion is made that the deformation parameter $\hbar$ in the series defining $a\star b$ for $ a,b \in \cA$ can be specialized to some set of values in the coefficient ring.   When this is possible for a particular UDF and algebra we will say that the UDF is effective on that algebra. An important example is that of the deformation of $\C[x,y]$ to the first Weyl algebra, discussed in \S\ref{Weyl}.

\begin{theorem}[\textsc The Basic UDF]\label{firstUDF}
 If $D_1, \,D_2$ are commuting derivations of an associative algebra $\cA$ over a ring $\bfk$ containing  the rationals, $\Q$, then
$ \exp \hbar(D_1\cpr D_2)$ integrates $(D_1\cpr D_2)$, i.e.,
  the multiplication defined on $\cA[[\hbar]]$ by
\begin{equation}\label{star product}
a\star b\, = \, \exp \hbar(D_1\cpr D_2)(\,a,b)\,   = \,\sum_{n=0}^{\infty}\frac{\hbar^n}{n!}D_1^na\cdot D_2^nb
\end{equation}
is associative and defines a formal deformation of  $\cA$.
\end{theorem}
\noindent \textsc{Proof.}
Setting 
 $F_i = D_1^i \cpr D_2^i$ in \eqref{basic}, what must be shown is that, for all $n$, one has the following relation among 3-cochains of $\cA$:
\begin{equation}\label{trivial}
\sum_{i=0}^n\frac{1}{i!(n-i)!}D_1^i\cpr D_2^i(D_1^{n-i} \cpr D_2^{n-i}) =
  \sum_{i=0}^n\frac{1}{i!(n-i)!}D_1^i(D_1^{n-i} \cpr D_2^{n-i})\cpr D_2^i .  
\end{equation}
By Leibniz' rule, the left  side of \eqref{trivial} can be written as 
\begin{multline}\label{left}
 \sum_{i=0}^n \frac{1}{i!(n-i)!} \sum_{j=0}^i \binom{i}{j} D_1^i \cpr D_2^jD_1^{n-i} \cpr D_2^{i-j}D_2^{n-i}\\
 =\sum_{i,j=0,\dots,n, \, i\ge j} \frac{1}{(n-i)!j!(i-j)!} D_1^i \cpr D_2^jD_1^{n-i} \cpr D_2^{n-j},
\end{multline}
and the right side of \eqref{trivial} can be written as
\begin{multline}\label{right}
 \sum_{i=0}^n \frac{1}{i!(n-i)!} \sum_{j=0}^i \binom{i}{j} D_1^{i-j}D_1^{n-i} \cpr D_1^jD_2^{n-i} \cpr D_2^i\\
 =\sum_{i,j=0,\dots,n, \, i\ge j} \frac{1}{(n-i)!j!(i-j)!} D_1^{n-j}\cpr D_1^jD_2^{n-i} \cpr D_2^i.
\end{multline}
On the right side of \eqref{right} one can  
now replace the dummy variable $i$ by $n-j$ and $j$ by $n-i$.  The sum continues to be over the same set of indices since $i\ge j$ if and only if $n-j \ge n-i.$ As $D_1$ and $D_2$ commute, the right side of \eqref{right} then becomes identical to the right side of \eqref{left}, proving the assertion. $\Box$
\smallskip

The Basic UDF has the following immediate extension.
\begin{cor}\label{sums}
If $D_i, \, i = 1, \dots,n$ are commuting derivations (not necessarily all distinct) of an algebra $\cA$ over $\Q$ and $c_{ij}, i,j = 1,\dots,n$ are central elements which are constants for all the $D_i$, i.e., central elements of $\cA$ with $D_kc_{ij}=0$ for all $i,j,k$ 
then
$\exp \sum_{i,j=1}^nc_{ij}( D_i \smile D_j)$ quantizes $\cA$. 
 \  $\Box$
\end{cor}

\noindent \textsc{Proof.} 
The case $n=2$ is immediate from the foregoing. Suppose it is true for a given $n$. A derivation which commutes with all $D_1, \dots , D_n$ remains a derivation in the deformed multiplication given by \eqref{star product}. Therefore, if we now have commuting derivations $D', D''$ which commute with all $D_, \dots, D_n$ and a central element $c$ which is a constant for $D_1, \dots, D_n, D', D''$, then $a \star' b :=
\exp c(D' \cpr_{\star} D'')(a,b)$  is a deformation of the previously defined star product, where $\,\cpr_{\star}$ indicates that the cup product is to be taken with the star multiplication.  As the derivations commute, the product of the exponentials involved is the exponential of the sum of the exponents. 
$\Box$ \smallskip

\section{Weyl and Clifford algebras}\label{Weyl}
Heisenberg, \cite{Heisenberg}, stated his uncertainty principle for a particle moving in one dimension  in the form
$$
pq-qp = \frac{h}{2\pi i},
$$ 
where  $p$ is the position of the particle and $q$ its momentum, considered as functions in phase space. Its algebraic meaning is codified in the first Weyl algebra, $\cW$, which can be defined as the tensor algebra $T(V)$ on a 2-dimensional vector space $V$ with basis $x,y$, modulo the ideal generated by $xy-yx-1$. This can be done more generally with a commutative unital ground ring $\bfk$ over the rationals. The resulting algebra will be denoted in the following by $\cW(\bfk)$.

\begin{theorem}
The Weyl algebra $\cW(\bfk)$ can be obtained as an effective deformation of the polynomial ring $\bfk[x,y]$,
by taking, in \eqref{star product},  $D_1= \px, \, D_2 = \py$. The infinitesimal deformation of $\bfk[x,y]$ is then $\px\cpr\py$, and the deformed star multiplication is given, for $a,b$ in $\bfk[x,y]$, by 
\begin{equation}\label{Weylstar}
a\star b = ab + \hbar \px a\py b + \frac{\hbar^2}{2!}\px^2a\,\py^2b + 
  \frac{\hbar^3}{3!}\px^3a\,\py^3b + \cdots .
 \end{equation} 
 \end{theorem}
 
 \noindent\textsc{Proof.} While the right side of \eqref{Weylstar} is formally an infinite series, it must terminate for any fixed $a$ and $b$. It follows that the deformation parameter, $\hbar$, can be specialized to any value in $\bfk$, so the deformation is effective.
 From \eqref{Weylstar}, one sees that $x\star y = xy + \hbar$ while $y\star x = yx = xy$, so the star commutator, $[x,y]_{\star}$, is $x\star y -y \star x = \hbar$. Setting $\hbar =1$, the Weyl algebra $\cW(\bfk)$ is thus a deformation of $\bfk[x,y]$. 
 $\Box$
 \smallskip
 
 The $n$th Weyl algebra, $\cW_n = \cW^{\otimes n}$,  is similarly an effective deformation of $\C[x_1, y_1, \dots, x_n, y_n]$ with infinitesimal $\sum_{i=1}^n \pl {x_i}\cpr\pl {y_i}$. It encodes the uncertainty principle for a particle moving in $n$ dimensions with position and momentum coordinates $p_i,q_i,\, i = 1, \dots, n$ and commutation relations $[p_i, p_j] = [q_i,q_j] = 0, [p_i,q_j] = \delta_{ij}$.
 
Clifford algebras can similarly be viewed as effective deformations of exterior algebras.  The latter are graded algebras, so one must follow Koszul's rule of signs.
In the smallest case, let $\cA$ be the four dimensional exterior algebra over a field $\bfk$ generated by a two dimensional vector space spanned by $x$ and $y$, each of which has degree 1. For simplicity, assume that the characteristic of $\bfk$ is not 2. Here $\px, \py$ are both derivations of $\cA$ of degree -1 which commute in the sense that $\px\py = -\py\px$, and each  has square equal to zero.  (The square of a derivation of odd degree always must be zero when the characteristic is not 2, since it commutes with itself.) The exponential $\exp(\px\cpr\py)$
reduces to $\id\cpr\id + \, \px\cpr\py$.  With the star product it induces one still has $x\star x  = y\star y = 0$ but now $x\star y = xy +1$, while $y\star x  = yx$, so $x\star y + y \star x = 1$. Therefore, $(x+y)^{\star2} =1$ and $ (x-y)^{\star2} = -1$, so the deformed algebra is the Clifford algebra
${\text Cl}_{1,1}(\bfk)$.

 \section{Fourier series}\label{Fourier}
Consider again the algebra   $C^{\infty}(\C^2)$ of smooth functions of two complex variables, $x,y$,  with commuting derivations $D_1=\px, D_2 = \py$ and infinitesimal deformation $\px \cpr \py$. It contains a subalgebra $\mathcal F
$ of for which the associativity of the resulting star product as well as its effectiveness  can be seen without the use of Theorem \ref{firstUDF}, namely that 
 generated by functions of the form $e^{\lambda x + \mu y}$.  
For these, applying \eqref{Weylstar} with $a = e^{\lambda_1 x + \mu_1 y},\,  b = e^{\lambda_2 x + \mu_2 y}$ one has first
$$
e^{\lambda_1 x + \mu_1 y}\star e^{\lambda_2 x + \mu_2 y} = 
e^{\mu_1y+\lambda_2 x}(e^{\lambda_1 x}\star e^{\mu_2 y})
$$
because $y$ is a constant for $\pl x$ and $x$ for $\pl y$, after which  computation gives
\begin{equation}\label{prestar}
e^{\lambda_1 x}\star e^{\mu_2 y} = e^{\hbar\lambda_1\mu_2}e^{\lambda_1 x +\mu_2 y}.
\end{equation}

Together one has
\begin{equation}\label{first star}
e^{\lambda_1 x + \mu_1 y}\star e^{\lambda_2 x + \mu_2 y} = 
e^{\hbar(\lambda_1\mu_2)}e^{(\lambda_1 + \lambda_2)x + (\mu_1 + \mu_2)y},
\end{equation}
where the right side is so written to exhibit its second factor as the undeformed multiplication and its first as the effect of the deformation.
\noindent From \eqref{first star} one may see first that it independently shows the associativity of this star multiplication  for functions of the special form $e^{\lambda x +\mu y}$ and second that the deformation is effective since one can specialize $\hbar$ to any complex value.  The star product may still be viewed as local, even though its value at a  point depends on the global behavior of the factors, since those factors are analytic and determined globally by their germs at any point. 

Functions of the form $e^{\lambda x + \mu y}$  are doubly periodic but those in $\mathcal F$, which are linear combinations of such, are generally only almost periodic because the periods of the summands may not be commensurable. Doubly periodic functions in $x$ and $y$ with fixed periods in each, say $2\pi$, can, with suitable restrictions,  be represented as double Fourier series using complex exponentials. Those with absolutely convergent double series form a subring, $\mathcal F_{\text{ac}}(2\pi,2\pi)$, on which the foregoing star product is still well-defined and effective: If in \eqref{prestar}, $\lambda_1$ and $\mu_2$ are both purely imaginary and and one replaces $\hbar$ by $i\hbar$ then the deformation factor will have absolute value 1, so the star product will still be absolutely convergent for all now real (previously purely imaginary) $\hbar$.  The functions in $\mathcal F_{\text{ac}}(2\pi, 2\pi)$ are generally not analytic and the deformed product is no longer local but is still effective. This should be extendable to Fourier transforms.

The foregoing allows a wide class of effective deformations stemming from the Basic UDF. While Kontsevich's UDF, \cite{Kontsevich:Poisson}, applies to any Poisson structure it remains purely formal; it can not define a star product on functions which are not analytic such as those above.

 \section{Jump deformations, rigidity, and stability}\label{jump}
The Weyl algebra is also an example of a jump deformation,  \cite[\S 7, p. 20ff]{GS:AlgebraicCohomology}, \cite{GG:qWeyl},  i.e., one in which the algebras defined for all values of $\hbar$ other than zero are  isomorphic. In a jump deformation, the infinitesimal of the deformation becomes a coboundary in the star multiplication, see e.g., \cite{GG:qWeyl}.   

The cohomology of a $\bfk$ algebra $\cA$ and that of a deformation, $\cA_{\hbar}$, of $\cA$ are closely related, \cite{G:Def4}, as follows.    The definition of the star product as an infinite series in $\hbar$ requires that the underlying coefficient ring of $\cA_{\hbar}$ be extended at least to the power series ring $\bfk[[\hbar]]$, but to remain over a field, let it be extended to Laurent series, $\bfk((\hbar))$, in which a finite number of negative powers of $\hbar$ are allowed.  An $n$-cocycle $z$ of $\cA$ is liftable to $\cA_{\hbar}$ if there are $n$-cochains $z_2, z_3, \dots$ such that $z + \hbar z_1 + \hbar^2z_2 + \cdots$ is an $n$-cocycle of $\cA_{\hbar}$.  The cohomology of $A_{\hbar}$  then consists of the liftable cocycles modulo those which lift to coboundaries. 

By the HKR theorem, the cohomology of $\C[x,y]$ is isomorphic, as a vector space, to the exterior algebra on $\px, \py$. Both $\px$ and $\py$ lift, in fact to themselves, as derivations of $\cW$, but become inner since, e.g., in $\cW$ one has $\py y = \hbar^{-1}[x,y]_{\star} = 1$, whence $\py a = \hbar^{-1}[x,a]_{\hbar}$ for all $a \in \cW$. It follows that the cohomology of $\cW$ vanishes in all positive dimensions, in particular, it is absolutely rigid in the terminology of \cite{G:Def}, or in the better terminology of Froelicher and Nijenhuis, \cite{FN}, stable.
This was first proven in Ramaiyengar Sridharan's thesis under S. Eilenberg, \cite{Srid:thesis}.  

As the Weyl algebra is a stable deformation of the polynomial ring in two variables, which itself is not stable, this might suggest that classical physical laws tend to deform to stable ones but the laws themselves do not change, it is our understanding of them that evolves. It is also remarkable that, at least in this context, uncertainty seems to be necessary for stability. 

If $\cA_1, \cA_2$ are algebras over the same field $\bfk$ with vanishing regular
 cohomology (i.e., with coefficients in themselves, cf \S \ref{History}) in all positive dimensions, then the same is true 
of their tensor product, a special case of \cite[Theorem 4, p. 207]{Mac Lane:Homology}. It follows that if in $\mathbb C^{\infty}[q_1,\dots,q_d,\,p_1,\dots, p_d]$ one takes as infinitesimal deformation $i\hbar\sum_{i=1}^d(\partial q_i \cpr \partial p_i)$, then the resulting deformed algebra, in which $q_ip_i -p_iq_i = i\hbar$ for all $i$ but $p_i$ and $q_j$ still commute for $i \ne j$, has no regular cohomology in positive dimensions. In particular, it is stable. 

By contrast, a tensor product of rigid algebras need not be rigid. As an example, consider the {\em twisting of a tensor product}: Suppose that we have algebras $\cA_1, \cA_2$ defined over a ring $\bfk$ containing $\Q$ with respective derivations $D_1, D_2$. 
Extend $D_1, D_2$  to derivations of $\cA_1 \otimes_{\bfk} \cA_2$ by setting $D_1(a_1 \otimes a_2) = D_1a_1\otimes a_2$ and $ D_2(a_1\otimes a_2) = a_1 \otimes D_2a_2$. These extensions commute, so $D_1 \cpr D_2$ can be exponentiated to a full deformation of $\cA_1 \otimes \cA_2$.  If either $D_1$ or $D_2$ is inner  then the infinitesimal deformation $D_1 \cpr D_2$ will be a coboundary and the resulting deformation will be trivial.   However, if neither is inner, then this is an example where $\cA_1 \otimes_{\bfk} \cA_2$ can be deformed even when both $\cA_1$ and $\cA_2$ are rigid.

Theorem \ref{firstUDF} was  first explicitly stated not only for associative algebras but for all composition complexes as defined in
\cite[Lemma 1, p.14]{G:Def3} but is already implicit in two important papers on quantum mechanics, Groenewold, \cite{Groenewold}, 1946, and Moyal, \cite{Moyal}, 1949. The algebra being deformed in both is the algebra of observables $\mathbb C [p,\,q]$ in phase space, where, as remarked, before deformation the position variable $q$ and momentum variable $p$ commute, but after deformation satisfy Heisenberg's uncertainty relation $qp-pq = i\hbar$. 

Infinitesimal deformations of $C^{\infty}(\mathbf R^d)$ of the form  $\sum_{i<j}c_{ij}f_i\pl i\cpr f_j\pl j$, with all $c_{ij} \in \R$  and each $f_i$ a smooth function only of the one variable $x_i$, will be called basic. From  the corollary to Theorem \ref{sums},  
the following is immediate.
\begin{theorem}\label{sums3} 
Basic infinitesimal deformations of $\Cinf$ are integrable by exponentiation.  \quad $\Box$
\end{theorem} 

\noindent
To shorten some lengthy expressions, in what follows we will write $(i|j)$ for $\pl i\cpr\pl j$,\, $(ij|k)$ for $\pl i\pl j\cpr\pl k$,\, $(i|j|k)$ for $\pl i\cpr\pl j\cpr\pl k$, and likewise for similar expressions.
With this notation, simple computation yields the following.

\begin{lemma}\label{prelim to F_1squared}For all $1 \le i < j \le d, 1 \le \, k < l \le d$, one has
\begin{align*} \label{products}
a_{ij}(i|j)a_{kl}(k|l)  &=  a_{ij}a_{kl}(ik|jl) \quad \text{if }
i\neq k, j \neq l\\  
a_{ij}(i|j)a_{il}(i|l) &= a_{ij}\pl ia_{il}(i|jl) + a_{ij}a_{il}(ii|jl) \quad \text {  if $j \neq l$;} \\
a_{ij}(i|j)a_{kj}(k|j)  &=  a_{ij}\pl ja_{kj}(ik|j) + a_{ij}a_{kj}(ik|jj)\quad\text{if $i \neq k$;} \\
a_{ii}(i|i)a_{ii}(i|i)  &=  \pl ia_{ii}\pl ia_{ii}(i|i) + a_{ii}\pl ia_{ii}(ii|i) + a_{ii}\pl ia_{ii}(i|ii) + a_{ii}a_{ii}(ii|ii). \Box 
\end{align*}
\end{lemma}

This lemma allows one to write all the powers, and hence the exponential, of an infinitesimal $F_1 = \sum_{ij}a_{ij}(i|j)$ in terms of the $a_{ij}$ and their various derivatives. In particular, for a basic infinitesimal $F_1 = \sum_{i<j}c_{ij}f_i\pl i\cpr f_j\pl j$, writing $a_{ij}$ for $c_{ij}f_if_j$, these relations allow one to write $F_1^2$, and hence all powers $F_1^n$ as well as $\exp{F_1}$,  in terms of the $a_{ij}$ and their derivatives without reference to the $f_i$.  However, if    $ F_1 = \sum_{1\le i < j \le d}  a_{ij}(i|j)$, is not basic then $\exp{F_1}$ will generally not be an integral of $F_1$. For it to be an integral, by \eqref{Basic2} we must have $F_1\circ F_2= \delta(-\frac{1}{2}F_1^2)$, but it is the difference between the left and right that gives rise to the obstruction equations discussed later. 
  
 \section{Differential subcomplex}
There is a differential subcomplex,  $C_{\text{diff}}^{\bullet}(\Cinf)$, of the Hochschild complex  $C^{\bullet}(\Cinf)$ of $\Cinf$ generated by the derivations $\pl 1, \dots \pl{d}$ with respect to its coordinates $x_1,\dots,x_d$; . Its 0-cochains are just the elements of $\Cinf$.
The module of n-multiderivations 
  is spanned by the $n$-fold cup products of these.   Simple 1-cochains are of the form $a\pl{i_1}\pl{i_2}\cdots \pl{i_r}$, where $a\in\Cinf$, $r$ is arbitrary, and there may be duplications amongst the indices.  The module of 1-cochains is composed of   sums of such cochains of various orders. Simple $n$-cochains are cup products of simple 1-cochains. 
The differential order of such a cup product  is the sum of the orders of its cup factors. A general $n$-cochain is a sum of simple ones, possibly of different differential orders; it is homogeneous if the orders are the same. 

The Hochschild coboundary operator preserves differential order. While the differential subcomplex is closed under the composition product, the composition product of homogeneous cochains is generally no longer homogeneous.  The HKR theorem implies that the inclusion of the differential subcomplex into the full Hochschild complex induces an isomorphism of cohomology.

Using the notation preceding Theorem \ref{sums3}, the HKR theorem also implies that $H^r(\cA)$ can be identified with the $\cA$ module spanned (using the notation of \S\ref{jump}) by all forms  $(i_1|i_2|\cdots |i_r)$  with $i_1< i_2 < \cdots < i_r$, since the class of any skew form also has a unique representative in that module. 

\section{Poisson structure}\label{Poissonstructure}
A Poisson structure on a commutative algebra $\cA$ is a second, 'bracket', multiplication
$[a,b], \, a,b \in \cA$, which is skew, $[a,b] = -[b,a]$, a biderivation, $[ab,c] = a[b,c] +b[a,c]$, the same holding on the right by skewness, and which satisfies the Jacobi identity,
\begin{equation}\label{Jacobi}
J(a,b,c):=[a, \, [b, \,c]] + [b, \, [c, \, a]] + [c, \, [a, \, b]] = 0,
\end{equation}
making $\cA$ with this second multiplication into a Lie algebra.  
The Poisson structure is an infinitesimal deformation of the algebra $\cA$ since a biderivation of a commutative algebra is clearly 2-cocycle.  
The `Jacobiator', $J$, 
is skew in all three variables.  Let $\frak g$ be a finite dimensional Lie algebra over a field $\bf k$ of arbitrary characteristic with basis $x_1, \dots, x_d$. The symmetric algebra on the underlying vector space of $\frak g$, which is just the polynomial ring $\bf k[x_1,\dots, x_d]$, carries a Poisson structure defined by extending the Lie bracket on $\frak g$ to be a biderivation of $\bf k[x_1, \dots , x_d]$. 

A Poisson manifold $\cM$ is one with a Poisson structure on its algebra of smooth functions. The dual $\frak  g^*$ of a finite-dimensional Lie algebra $\frak g$ over $\R$ or $\C$ is an example. It carries a Poisson bracket called variously a Kirillov-Poisson, Lie-Poisson,  or Kirillov-Kostant-Souriau (KKS)  structure described as follows. If $f$ is a differentiable function on a finite dimensional vector space $V$ and $p$ a point of $V$, then the differential of $f$ evaluated at $p$, $d_pf$, is an element of the dual space $V^*$. 
 Suppose now that $f, g$ are functions on the dual, $\frak g^*$, of a finite dimensional Lie algebra $\frak g$ and $\xi$ a point of $\frak g^*$.  One can then take the Lie product $[d_{\xi}f,\,d_{\xi}g]_{\frak g}$ in $\frak g$ of their differentials at $\xi$ and define the value at $\xi$ of the KKS Lie product $[f,g]_{\frak g^*}$ of $f$ and $g$ to be $\xi([d_{\xi}f,\,d_{\xi}g]_{\frak g})$. 
To calculate the infinitesimal deformation which the Kirillov Poisson structure defines on $\frak g^*$, choose a basis $x_1,\dots, x_d$ of $\frak g$, let  $\xi_1,\dots, \xi_d$ be its dual basis in $\frak g^*$, let $\xi  = \sum a_i\xi_i$ be point in $\frak g^*$ and suppose that $[x_i,x_j] = \sum c_{ij}^kx_k$. Then 
\begin{equation*}
df = \sum_1^dx_i\frac{\pl f}{\pl {\xi_i}},\quad
[df, dg] = \sum c_{ij}^kx_k\frac{\pl f}{\pl {\xi_i}}\frac{\pl g}{\pl {\xi_j}}, \, \text{and}\,\, \xi([df,dg]) = \sum c_{ij}^ka_k\frac{\pl f}{\pl {\xi_i}}|_{\xi}\frac{\pl g}{\pl {\xi_j}}|_{\xi}
\end{equation*}
where all indices run from 1 to d and the partial derivatives in the last term on the right are evaluated at $\xi = \sum a_i\xi_i$.  To represent this in the form 
$$
(\sum \lambda_{rs} (x_r \cpr x_s)(f,g))|_{\xi} =
\sum \lambda_{rs} x_r(f)|_{\xi}x_s(g)|_{\xi}
$$
take $f= \xi_i, \, g= \xi_j$; comparing with the preceding gives
\begin{equation}\label{bivector} 
\lambda_{ij} = \sum c_{ij}^ka_k.
\end{equation}
The KKS product of linear functions on $\frak g^*$ is linear. Conversely, any Poisson manifold  in which the product of linear functions is linear arises in this way.
If $\xi_i$ is a central element of $\frak g$ then all $c_{ij}^k$ are zero so 
 if $\frak g$ is Abelian, then the KKS structure on $\C^{\infty}(\frak g*)$ vanishes and induces no deformation.  Nevertheless, as the elements $x$ of the Lie algebra $\frak g$ act as derivations on $\C^{\infty}(\frak g^*)$ by $\pl {x}f|_{\xi}=
 x(f)_{\xi} = \xi[x, d_{\xi}f]_{ g}$,  if $\frak g$ has two elements which commute then one can apply the Basic UDF to deform $\C^{\infty}(\frak g^*)$. 
 
 For any $\xi \in \frak g^*$ and differentiable $f$  on $\frak g^*$ there is a linear function with the same differential at $\xi$ as $f$, so the Poisson bracket on $\C^{\infty}(\frak g^*)$ is already determined by that on its linear functions.  
 
\section{Coll-Gerstenhaber-Giaquinto UDF}\label{non-com}
Theorem \ref{sums3} provides another approach to the Coll-Gerstenhaber-Giaquinto (CGG) UDF, \cite{CGG:noncomm}, the first with non-commuting derivations.  In the following, $[D_2]_n$ will denote the ``descending factorial", $[D_2]_n = D_2(D_2-1)\cdots(D_2-n+1)$ and $n$ will mean $n \id_{\cA}$, so, e.g., $[D_2]_2$ means $D_2(D_2- \id_{\cA})$.
\begin{theorem}\label{CGG}
Suppose that $D_1,D_2$ are derivations of an algebra $\cA$ over $\Q$ such that $[D_1,D_2] =D_1$. Then  $D_1\cpr D_2$ is an integrable  infinitesimal deformation of $\cA$ with integral 
\begin{multline*}
e(\hbar, D_1,D_2): = \\id_{\cA}\cpr id_{\cA}+ \hbar D_1\cpr D_2 +(\hbar^2/2!)D_1^2\cpr [D_2]_2 + \cdots + (\hbar^n/n!)D_1^n\cpr[D_2]_n + \cdots \quad .
\end{multline*}
\end{theorem}
\noindent \textsc{Proof.}
The UDF asserted here is essentially a sequence of formal identities in $D_1,D_2$  within the universal enveloping algebra of the two-dimensional Lie algebra generated by $D_1,D_2$. 
If we can exhibit an explicit associative algebra $\cA$ with derivations $D_1,D_2$ having an isomorphic universal enveloping algebra then this will serve as a model in the sense that any proposition about the universal enveloping algebra which holds in this model must be true in general.
 We can model $D_1$ and $D_2$ as derivations of the algebra of smooth functions of $x$ and $y$ by mapping $D_1$ to  $e^{-y}\pl x$ and  $D_2$ to $\pl y$, but must show that the natural map of the universal enveloping algebra of the Lie algebra generated by $D_1$ and $D_2$ onto the algebra of operators generated by  $e^{-y}\pl x$ and  $\pl y$ is an isomorphism. The commutation relation between $D_1$ and $D_2$ implies that the universal enveloping algebra is, as a module over the ground ring, free with generators all monomials of the form $D_1^mD_2^n$, with $m,n \ge 0$. The analogue clearly also holds for the algebra of operators generated by $e^{-y}\pl x, \pl y$, so they are isomorphic modules.
 Now $D_1\cpr D_2 = e^{-y}\pl x \cpr \pl y$, but this can also be written as   $\pl x \cpr e^{-y}\pl y$, which is a basic infinitesimal with integral 
\begin{equation}\label{CGG2}
\operatorname{exp}(\hbar\pl x \cpr e^{-y}\pl y) = \sum_{n=0}^{\infty}(\hbar^n/n!)(\pl x)^n \cpr (e^{-y}\pl y)^n.
\end{equation}
It is a simple induction to show that
for every $\lambda \in \mathbb{C}$ and positive integer $n$ one has
$$(e^{-\la y}\pl y)^n = e^{-ny}\pl y(\pl y-\la)(\pl y-2\la)\cdots (\pl y-(n-1)\la).
$$
In particular, for $\la = 1$, one has
$(e^{-y}\pl y)^n = e^{-ny}[\pl y]_n.$  The right side of \eqref{CGG2}  can therefore be written as \,
 $\sum_{n=0}^{\infty}(\hbar^n/n!)(e^{-y}\pl x)^n \cpr [\pl y]_n$,\, proving the theorem.
 $\Box$\smallskip
 
 An important property of this UDF is  that it gives an explicit integral for the KKS structure (\S\ref{Poissonstructure}) on the dual of the unique non-abelian solvable two-dimensional Lie algebra $\mathfrak s$.  One can take a basis, $\{x_1, x_2\}$, for $\mathfrak s$ with $[x_1, x_2] = x_1$. With this, in \eqref{bivector}, one has $c_{12}^1 = 1, c_{12}^2 =0$, which completely determines $\Pi$, since it is skew; $\Pi = \pl {x_1} \cpr \pl {x_2} - \pl {x_2} \cpr \pl {x_1}$.  Theorem \ref{CGG} integrates this since it is cohomologous to $ 2\,\pl {x_1} \cpr \pl {x_2}$. The integral is rational and computable term by term.

\section{The primary obstruction}\label{primeobs}
As remarked in \S \ref{Poissonstructure}, a biderivation of a commutative algebra is a 2-cocycle, hence an infinitesimal deformation. As a partial converse, one has the following.
\begin{lemma}\label{3delta}
Let $f$ be a skew 2-cocycle of a commutative algebra $\cA$ over a ring $\bfk$ in which 2 is a unit. Then $f$ is a biderivation.
\end{lemma}
\noindent \textsc{Proof.}
For a skew 2-cochain $f$ of a commutative algebra one has the identity
\begin{equation*}
(\delta f)(a,b,c) +(\delta f)(c,a,b) - (\delta f)(a,c,b) = 
2[af(b,c) - f(ab,c)  + bf(a,c)].\, \Box
\end{equation*}
\smallskip
 
\noindent Lemma \ref{3delta} suggests that an alternating  $n$-cocycle (one which changes sign under an odd permutation of its arguments) of a commutative algebra  should be a multiderivation when $n!$ is a unit. 
Recall from the HKR theorem that every element of $H^2(\Cinf)$  has a unique skew representative 2-cocycle, which by Lemma \ref{3delta} is a biderivation, and from \S \ref{BasicUDF} that the primary obstruction to an infinitesimal deformation $F_1$ of an algebra is $F_1\circ F_1$.

In what follows, for a function $f$ of three variables, set
$$ \sum_\circlearrowleft f(a,b,c) := f(a,b,,c)+f(b,c,a)+f(c,a,b).$$
Note that if $f$ is skew in any pair of variables then  $ \sum_\circlearrowleft f$ is skew in every pair.

\begin{theorem}\label{Poisson}
1. Let $\Pi$ be a skew 2-cocycle of a commutative algebra $\cA$ defined over a ring $\bfk$ in which $3!$ is a unit. Viewing $\Pi$ as an infinitesimal deformation of $\cA$, if its primary obstruction, $\Pi \circ \Pi$, vanishes then $\Pi$ is Poisson. 2. An infinitesimal deformation of a commutative algebra $\cA$ defined over $\Q$ is integrable to a full deformation if and only if its primary obstruction vanishes.
\end{theorem}

\noindent \textsc{Proof.} 
Assertion 2. follows immediately from 1.  since if $\Pi$ is Poisson then by \cite{Kontsevich:Poisson} it is integrable, so in particular its primary obstruction vanishes. Suppose now that the primary obstruction to $\Pi\circ\Pi$ vanishes, i.e., that it is a coboundary.
In view of Lemma \ref{3delta} it is sufficient now to show that   $\sum_\circlearrowleft\Pi\circ\Pi = 0$. Since $\Pi\circ\Pi$ is a coboundary so is $\sum_\circlearrowleft\Pi\circ\Pi$. Since $\Pi$ is skew one has that $\Pi\circ\Pi(a,b,c)=\Pi(\Pi(a,b),c)- \Pi(a, \Pi(b,c))$ is skew in it first and third variables, $a$ and $c$. It follows that $\sum_\circlearrowleft\Pi\circ\Pi$, which is again a coboundary, is skew in all three variables. From the decomposition of commutative algebra cohomology introduced in \cite{GS:Hodge} it must be identically zero: There it was shown that if $n!$ is a unit then the $n$th cohomology $H^n$ of a commutative algebra with coefficients in a commutative module decomposes into a direct sum
$$H^n = H^{1,n-1} + H^{2,n-2} + \cdots + H^{n,0}$$
where the `top' component, $H^{n,0}$ consists of skew multiderivations each the unique representative of its cohomology class; a skew mutiderivation therefore can not be a coboundary unless it is identically zero. $\Box$

 \smallskip
  
 The UDF given by Kontsevitch in \cite{Kontsevich:Poisson}
has the property that 
  the higher order terms come from various partial derivatives of the coefficients in the infinitesimal deformation, so the full deformation it constructs  will be trivial in any open set where all coefficients of its infinitesimal vanish. 

\section{Moduli spaces and modular groups}\label{modular}
A deformation of an algebra or other structure is usually presented initially as dependent on one or more parameters. 
 The moduli space of the deformation is its parameter space with points representing isomorphic structures identified. It carries the quotient topology, in which a set is open if and only if its preimage is open. The parameter space can sometimes be discrete, as  in \cite{G:PathAlgebras}, and the moduli space can even be non-Hausdorff. This is the case with a jump deformation, where the moduli space consists of two points, one closed representing the original algebra, and one open representing the deformed one. A deformation problem will have an associated modular group when specializations of the deformation parameter yielding structures isomorphic to the original are discrete.  The modular group then operates on the parameter space to identify points corresponding to isomorphic objects.

An example of a modular group associated with the deformation of the polynomial ring in two variables to the quantum plane is given below.  It is analogous to the classical modular group  associated to the deformation of compact Riemann surfaces of genus one or complex tori. For these and compact Riemann surfaces of higher genera there is a Teichm{\"u}ller space whose points are essentially the conformal equivalence classes of the surfaces each with a choice of integral generators for the homology group of its underlying manifold. In the case of genus one, choosing one surface with marked generators of its integral homology groups, the others are then marked by the elements of the special linear group $\mathbf{SL}(2,\Z)$ which is isomorphic to the free product $\Z/4\star \Z/6$. 
 In detail,
 every complex torus can be represented as the quotient $\mathbb C/L$  of the complex plane $\C$ by a lattice $L$, which one may view as  
 generated by an unordered pair of real non-collinear real vectors 
 $\begin{bmatrix}x_1\\y_1\end{bmatrix}$,\,
 $\begin{bmatrix}x_2\\y_2\end{bmatrix};$
 they can be given a natural order by requiring that the determinant of 
 $\begin{bmatrix}
 x_1 &  x_2\\
 y_1 &  y_2\\
 \end{bmatrix}$
 be positive. The angle from the first vector to the second in the positive direction is then less than $\pi$.
The parallelogram they generate is a fundamental domain for $L$ considered as a group operating on $\C$ and the torus may be viewed as the fundamental domain with opposite sides identified.
The sides become a pair of circles intersecting in a single point and are generators for the first homology group of the torus. The special linear group $\mathbf{SL}_2(\Z)$ operates on these pairs of vectors, identifying pairs which generate isomorphic latices. Each of the lattices has an automorphism of order two obtained by rotating the parallelogram through a half circle around its center. Two have larger groups of automorphisms, namely that in which the parallelogram is a square where the automorphism group is $\Z/4$ and that in which it has equal sides separated by an angle of $\pi/3$ whose automorphism group is $\Z/6$. Accordingly, $\mathbf{SL}_2(\Z)$ is an amalgam; it is the free product of $\Z/4$ and $\Z/6$ with the subgroups of order 2 identified, $\Z/2 \star_{\Z/2}\Z/6$, cf \cite{G:Discontinuous}. Its center is the consolidated group $\Z/2$. As matrices in $\mathbf{SL}_2(\Z)$ it consists of the identity and its negative.

Viewing vectors in the plane as complex numbers, any complex torus $\C/L$ is biholomorphic to one where the first generating vector of $L$ is 1 and the second therefore a complex number $\tau$ with positive imaginary part. The parameter space for complex tori is thus the upper half plane. The projective special linear group $\mathbf{PSL}_2(\Z)$ consisting of all maps from $\C$ to itself of the form $\tau\mapsto (a\tau+b)/(c\tau+d),\, a,b,c,d \in Z, ad-bc =1$ acts on it to identify values of $\tau$ which represent the same complex torus. One can visualize the deformation of a torus as a stretching of the fundamental domain onto a new one by the motion of $\tau$ but the underlying structure of a manifold is preserved. This classical modular group, associated to the deformation of complex tori,                                                                                                                                                                 is the quotient of $\mathbf{SL}_2(\Z)$ by its center and is isomorphic to $\Z/2 \star \Z/3$.  The classical moduli space is smooth and biholomorphic to the punctured sphere, i.e., Riemann sphere with one point removed.

 For an analogous purely algebraic modular group and moduli space consider  now the deformation of $\mathbb C[x,y]$ with infinitesimal $x\pl x \cpr y\pl y$ and star product given by $\exp (\hbar \,x\pl x \cpr y\pl y)$.  This is again an example of the twisting of a tensor product of $\bfk [x]$ with itself. Previously the derivation chosen in each tensor factor was $\pl x$ while it is now $x\pl x$. 
With the resulting 2-cocycle $x\pl x \cpr y\pl y$ one has
\begin{equation*}
x\star y = xy  +\hbar (x\px)x\,(y\py)y+ \frac{\hbar^2}{2}(x\px)^2x\,(y\py)^2 y + 
\frac{\hbar^3}{3!}(x\px)^3x\,(y\py)^3y+ \cdots.
\end{equation*}
As $(x\px)^nx =x$ for all $n\ge 0$, and similarly with $y$, one has 
 $x\star y = e^{\hbar}xy$, while $y\star x =yx = xy$.  Therefore, with this deformation, 
 \begin{equation}\label{exphbar}
x\star y = e^{\hbar}y \star x.
 \end{equation}
   Writing $q$ for $e^{\hbar}$ one has $x\star y = qy\star x$, the defining equation of the quantum plane. The parameter space is all of $\mathbb C$.  The modular group, $\mathbb G$, is generated by  transformations $\sigma:\hbar\mapsto \hbar+2n\pi i, n \in \Z$ and $\tau:\hbar \mapsto -\hbar$,
the latter
corresponding to the 
interchange of $x$ and $y$.
There is an exact sequence 
\begin{equation}\label{exact}
0\longrightarrow \Z \longrightarrow \mathbb G \longrightarrow \Z/2 \longrightarrow 0
\end{equation}
which is split by sending the non-identity element of $\Z/2$ to $\tau$, so $\mathbb G$ is the semidirect product $\Z/2 \rtimes \Z$, the infinite dihedral group frequently denoted $\mathbf D_{\infty}$. There is an underlying structure preserved by the deformations, namely that of the tensor product of two vector spaces. Two algebras amongst the foregoing deformations have non-trivial automorphisms which preserve this structure, namely the original polynomial algebra and the anticommutative one with $xy =-yx$. In each case the automorphism is the interchange of tensor factors, so by \cite{G:Discontinuous} the modular group for this deformation problem is $\Z/2 \star \Z/2$ which is identical with $\mathbf D_{\infty}$.
  
The quotient of $\C$  by $\Z$ (appearing in \eqref{exact} as multiples of $\sigma: \hbar \to \hbar + 2\pi i$) is an infinite cylinder, which is conformal by the exponential map to $\C\backslash 0$, i.e., $\C$ with zero removed. This is  evident in \eqref{exphbar}.   The group $\Z/2$ operates on $\C\backslash 0$ by sending $z$ to $1/z$, in \eqref{exphbar}, $e^{\hbar} \to e^{-\hbar}$, so the quotient, the moduli space, is again the punctured sphere.
 This moduli space, like the classical one, has no singularities. However, the moduli space for Riemann surfaces of genus 2, which has complex dimension three, has a singularity at the point representing the surface given by $y^2 = x^6 - 1$; there is a discussion of the modular group in genus 2 in \cite{Schiller}. 
 
 The deformations of $\C[x,y]$ to the Weyl algebra, to the quantum plane, $W_{qp} = \C\{x,y\}/(xy-qyx)$ and to the $q$-Weyl algebra $W_q = \C\{x,y\}/(xy-qyx - 1)$ are discussed in \cite{GG:qWeyl} where the cohomology of each of the latter two algebras are computed. 
 The generic $q$-Weyl algebra is not rigid while the Weyl algebra, its specialization at $q=1$ is absolutely rigid or stable. Note, as observed in \cite{GG:qWeyl}, that the algebra $W_q(\hbar) =  \C{x,y}/(xy-qyx - \hbar)$ is isomorphic to $W_q$: replace $\hbar$  by $x^n\hbar$ for any $n >0$ and $y$ by $y+x^{n-1}\hbar/(1-q)$.  The results of \cite{GG:qWeyl} suggest that the moduli space of the quantum Weyl algebras consists of two copies of that of the quantum plane, the first with its usual topology and the second with the discrete topology, but the points of the first are not closed in the whole space, the closure of a point in the first copy containing the corresponding point in the second.                                                                                                                                                                                                                                                                                                                                                                                                                                                       

\section{Obstruction equations}
There are $\binom d3$ non-linear partial differential `obstruction equations'  in the coefficients $a_{ij}$ of an infinitesimal deformation $F_1 = \sum_{1\le i< j \le d}a_{ij}(i|j)$.
These arise from the difference between $F_1\circ F_1$ and $\delta(-\frac12 F_1^2)$; the primary obstruction to $F_1$ vanishes if and only they are satisfied. 
The terms in $F_1\circ F_1$ are in natural one-to-one correspondence with those in $F_1^2$ with $a_{ij}(i|j)\circ a_{kl}(k|l)$ in the former corresponding, in the notation of  Section \ref{BasicUDF}, to $a_{ij}(i|j)a_{kl}(k|l)$ in the latter.   
The terms of differential order 4 in $F_1\circ F_1$ and  $\delta(-\frac12 F_1^2)$ must coincide since that would be the case if all the coefficients were constants. While $F_1^2$ contains terms of differential order 2, these are cocycles and will not contribute to its coboundary. The obstruction equations arise, therefore, from the difference between the terms of order 3 of $F_1\circ F_1$ and $\delta(-\frac12 F_1^2)$.    
Analogous to the products in Theorem \ref{prelim to F_1squared}, there are the following four kinds of terms in 
$F_1\circ F_1:$
\begin{align*}
 (1)& \sum_{1\le i<j\le d, 1\le k<l\le d}a_{ij}(i|j)\circ a_{kl}(k|l), \, \text{where} \, i,j,k,l \, \text{are all distinct}\\
 (2)& \sum_{1\le i<j\le d, 1\le i<l\le d}a_{i|j}(i|j)\circ a_{il}(i|l), \, \text{where} \, j,l \, \text{are distinct}\\
 (3)& \sum_{1\le i<j\le d, 1\le k<j\le d}a_{ij}(i|j)\circ a_{kj}(k|j), \, \text{where} \, i,k \, \text{are distinct}\\
 (4)& \,\,\,\,\,\,\,\,\,\,\,\sum_{1\le i<j\le d}\,\,\,\,\,\,\,\,\,\,\,a_{ij}(i|j)\circ a_{ij}(i|j).
 \end{align*}
Each composition product above gives rise to two terms of the form $b_{ijk}(i|j|k)$, where there may be a single repetition among $i,j,k$. A brief examination will show the following.
\begin{lemma}
 The sum of all the terms $b_{ijk}(i|j|k)$ in $F_1\circ F_1$ in which there is a repetition amongst $i,j,k$ is precisely $\delta(-\frac12 F_1^2)$. 
 $\Box$
 \end{lemma}
 It follows that $F_1\circ F_1 -\delta(-\frac12 F_1^2)$ is just the sum of those terms without repetitions.
The vanishing of the primary obstruction to a general infinitesimal $F_1$ is therefore equivalent to having the sum of those terms without repetitions be a coboundary.  
When the dimension $d = 2$, the HKR theorem implies that there  can be no obstructions; every infinitesimal deformation is then integrable. For 
every obstruction must lie in $H^3(C^{\infty}(\R^2))$, whose elements are uniquely represented by cocycles of the form $\sum_{1\le i <j<k \le d}a_{ijk}(i|j|k)$, but there can be no such cocycles when $d=2$. When $d \ge 3$, for every triple $(i,j,k)$ with $1\le i <j <k \le d$ a term of the form $b_{ijk}(i|j|k)$ appears amongst the composition products above, e.g., in $a_{ii}(i|i)\circ a_{jk}(j|k)$.  If $(i',j',k')$ is a permutation of $(i,j,k$) then $(i'|j'|k')$ is cohomologous to $(i|j|k)$ if the permutation is even, and to $-(i|j|k)$ if the permutation is odd. Therefore, the sum of the terms without repetitions is cohomologous to a unique sum of the form $\sum_{1\le i <j<k \le d}a_{ijk}(i|j|k)$. Here the $a_{ijk}$ are expressions in the coefficients $a_{ij}$ of $F_1 = \sum_{1\le i< j \le d}a_{ij}(i|j)$ and their partial derivatives with respect to the variables $x_1,\dots,x_d$. 
 For the primary obstruction to vanish, each of the $\binom d3$ coefficients $a_{ijk}$ must vanish, giving rise to $\binom d3$ non-linear partial differential equations in the coefficients of $F_1$. 

\section{The obstruction equation in dimension $d=3$}\label{obstruction}
In dimension 3, an infinitesimal deformation has the form $F_1= a_{12}(1|2) + a_{13}(1|3) + a_{23}(2|3)$ and  there is just one obstruction equation. It has a remarkable property, expressed in Theorem \ref{varphi}, below. 
One has 
\begin{align*}\label{FF}
F_1\circ F_1 = &\\
{}& a_{12}(1|2)\circ a_{12}(1|2)+ a_{13}(1|3)\circ a_{13}(1|3) + 
a_{23}(2|3)\circ a_{23}(2|3) + \\
& a_{12}(1|2)\circ a_{13}(1|3) + a_{13}(1|3)\circ a_{12}(1|2) + \\
& a_{12}(1|2)\circ a_{23}(2|3) +	a_{23}(2|3)\circ a_{12}(1|2) + \\ 
& a_{13}(1|3)\circ a_{23}(2|3) +	a_{23}(2|3)\circ(a_{13}(1|3).
\end{align*}
After expanding the composition products above gives
\begin{align*}
F_1\circ F_1 - \delta(-\frac12 F_1^2) = &\\ 
&\,\,\,\,\,\,\,\,(a_{12}\pl1a_{13})(1|3|2) + (a_{13}\pl1a_{12})(1|2|3) +\\
&-(a_{12}\pl2a_{23})(1|2|3) + (a_{23}\pl2 a_{12})(1|2|3)+\\
&-(a_{13}\pl3a_{23})(1|2|3) - (a_{23}\pl3 a_{13})(2|1|3).
\end{align*}
Since
$(1|3|2) =  -(1|2|3) + \delta(1|23)$
and
$(2|1|3) =  - (1|2|3) - \delta(12|3),$ 
we can write the sum of the foregoing terms as 
\begin{multline*}
[-a_{12}\pl1a_{13}+a_{13}\pl1a_{12}-a_{12}\pl2a_{23} + a_{23}\pl2a_{12} -a_{13}\pl3a_{23} + a_{23}\pl3a_{13}](1|2|3)\\
+ \delta(a_{12}\pl1a_{13} + a_{23}\pl3a_{13})(1|23).
\end{multline*}

The single obstruction equation in dimension 3 is therefore
\begin{equation}\label{obs3}
-a_{12}\pl1a_{13}+a_{13}\pl1a_{12}-a_{12}\pl2a_{23} + a_{23}\pl2a_{12} -a_{13}\pl3a_{23} + a_{23}\pl3a_{13} = 0.
\end{equation}

In any region where none of $a_{12}, a_{13}. a_{23}$ vanish, \eqref{obs3} can be rewritten as
\begin{equation}\label{obs3'}
a_{12}^2\pl1(\frac{a_{13}}{a_{12}}) + a_{23}^2\pl 2(\frac{a_{12}}{a_{23}}) + a_{13}^2\pl 3(\frac{a_{13}}{a_{23}}) = 0.                                                                                                       
\end{equation}

\begin{theorem}\label{varphi}
If $F_1$ is an infinitesimal deformation of $C^{\infty}(\R^3)$ whose primary obstruction vanishes then the same is true of $\varphi F_1$ for any $\varphi$ in $C^{\infty}(\R^3)$. In particular, if a skew biderivation  $F_1$ of $C^{\infty}(\R^3)$ is Poisson, then so is $\varphi F_1$. 
\end{theorem}
\noindent \textsc{Proof.}  
In \eqref{obs3}, if $F_1$ is replaced by $\varphi F_1$, where $\varphi$ is a smooth function of $x_1,x_2,x_3$, then the terms in which $\varphi$ is differentiated cancel, so the summand on the left side is just multiplied by $\varphi^2$, which may be more easily seen in \eqref{obs3'}. Therefore, if it is satisfied by a triple of smooth functions $(a_{12}, a_{13},  a_{23})$ of $x_1, x_2, x_3$,  then $(\varphi a_{12}, \varphi a_{13}, \varphi a_{23})$ is also a solution.
$\Box$
\smallskip

The extraordinary property of dimension 3 presented in Theorem \ref{varphi} does not appear in dimension 4 and therefore neither in any higher dimension. In dimension $d=4$, to see whether an infinitesimal  deformation 
$F_1 = \sum_{1\le i < j \le4} a_{ij}(i|j)$ has vanishing primary obstruction, we need only examine those terms in $F_1\circ F_1$ which have differential order three. The sum of these can be written as $b_{123}(1|2|3) + b_{124}(1|2|4) + b_{134}(1|3|4) + b_{234}(2|3|4) + $(coboundaries). The four obstruction equations assert the vanishing of the four $b_{ijk}$. It will be sufficient to show for one of them, say $b_{123}$, that vanishing for some $F_1$ does not imply that it also does so for all $\varphi F_1$, where $\varphi$ is a smooth function of $x_1, \dots, x_4$. The contributions to $b_{123}$ come first from those terms in $F_1\circ F_1$ of the form $a_{ij}(i|j)\circ a_{kl}(k|l)$ in which none of the indices $i,j,k,l$ is 4, and second from from $a_{14}(1|4)\circ a_{23}(2|3), a_{24}(2|4)\circ a_{13}(1|3)$, and
$a_{34}(3|4)\circ a_{12}(1|2)$.  In the first terms, if $F_1$ is replaced by $\varphi F_1$ then, as we have seen with the single obstruction equation  in dimension three, the contribution from those terms is just multiplied by $\varphi^2$. This is not the case with the second terms, so the statement fails  in dimension 4 and hence in all higher dimensions.  

\section{Poisson structures in dimension three}\label{Poisson3}
It is an open problem  to determine all Poisson structures on $\Cinf$ or equivalently to determine those infinitesimal deformations whose primary obstructions vanish, since by Theorem \ref{Poisson}  their skew forms are precisely the Poisson structures.  One can, however, construct a large set of such in the special case of dimension $d=3$. 
As before, letting $F_1 =  a_{12}(1|2) + a_{13}(1|3) +  a_{23}(23)$ it will be seen that any one of the functions $a_{12}, a_{13}, a_{23}$ can be chosen to be an arbitrary smooth function of $x_1, x_2, x_2$. Because of the symmetry in the variables suppose it is $a_{12}$.  By Theorem \ref{varphi} it will be sufficient then to give examples where $a_{12}=1$ but note that those so obtained need not be all. 
 Writing again $\pl i$ for $\pl {x_i}$,  \eqref{obs3} then becomes
\begin{equation}\label{obs3.1}
-\pl1a_{13} - \pl2a_{23}  - a_{13}\pl3a_{23} + a_{23}\pl3a_{13} = 0,
\end{equation}
which will be called the reduced obstruction equation in dimension 3. It can be rewritten as
\begin{equation}\label{obs3.2}
(\pl2 + a_{13}\pl3)a_{23} - (\pl3a_{13})a_{23} = -\pl1a_{13}.
\end{equation}
To find solutions to this equation, observe that if we first choose for $a_{13}$ an arbitrary smooth function of $x_1, x_2, x_3$ then \eqref{obs3.2} would  resemble a first order ordinary differential equation in a single variable $x$ of the form 
\begin{equation}\label{ODE}
\frac {dy}{dx} + p(x)y = q(x),
\end{equation}
with the differential operator $(\pl2 + a_{13}\pl3)$ in place of $\frac{d}{dx}$, $-(\pl3a_{13})$ in place of $p$, $-\pl1a_{13}$ in place of $q$, and where $a_{23}$ is in place of the unknown $y$.
Letting $P$ be any function with $ P' =p$,
 equation \eqref{ODE} has, as a particular solution,
\begin{equation}\label{part_soln}
y = \exp(-P(x))\int\exp (P(x))\,q(x)dx
\end{equation}
from which all others can be obtained by adding a solution to the homogeneous  equation $\frac {dy}{dx} + p(x)y = 0$
These all have the form
$c\,\exp(-P(x))$
where $c$ is a constant, so the general solution is
\begin{equation}\label{gen_soln}
y = \exp(-P(x))\int\exp (P(x))\,q(x)dx + c\,\exp(-P(x)).
\end{equation}
The only properties of the exponential function needed for \eqref{gen_soln} are that for any differentiable function $f$ of $x$ we have $(\exp f)' =(\exp f)f'$ and $exp(-f) = (exp f)^{-1}$. This determines the exponential function up to a sign; note that \eqref{gen_soln} would remain correct if the exponential were replaced by its negative.                                            
More generally, suppose that we have a commutative unital ring  $\mathcal R$ of arbitrary characteristic with  a derivation $\mathcal  D$ 
 and a function $\exp: \mathcal R \to \mathcal R$ such that (i) for any  $r \in \mathcal R$ we have $\exp (-r)= (\exp r)^{-1}$ and (ii) $\mathcal D (\exp(r)) = \exp r\,\mathcal D r$. If $p,q \in \mathcal R$ then the solutions to the equation $\mathcal D y +py = q$ have exactly the same form as in \eqref{gen_soln} but with the usual exponential function replaced by this $\exp$ and with the integral replaced by any fixed preimage under $\mathcal D$. One then has 
 \begin{theorem}\label{genODE}
 Let $\cR$ be a commutative unital ring with a derivation $\cD$ relative to which there is an `exponential function' $\exp$ such that for any  $r \in \mathcal R$ we have (i) $\exp  (-r) = (\exp r)^{-1}$ and (ii) $\mathcal D (\exp  r) = \exp(r)\,\mathcal D r$. If $P, q \in  \mathcal R$ and $\mathcal D P = p$ then a particular solution to
 \begin{equation}
 \mathcal D y +py = q
 \end{equation} 
 is given by
\begin{equation}
y = \exp(-P)\int(\exp P)\,q,
\end{equation}
from which all other solutions may be obtained by adding a function of the form $c\,\exp(-P)$, where $c$  is a constant for, i.e. annihilated by, $\cD$. $\Box$
\end{theorem}

Theorem \ref{genODE} may apply when instead of functions of one real variable we have functions of several. For `$\exp$' we may then take the usual exponential function. In prime characteristics one can use the Artin-Hasse exponential, for a discussion of which cf e.g. \cite{Matt:Exp}.
In the classical case of a single variable the derivation $\frac{d}{dx}$ is surjective on $\C^{\infty}(\R)$ so there is no restriction on $p$. This is not the case in \eqref{obs3.2} since it is required that $\pl 3a_{13} = (\pl 2 + a_{13}\pl 3)P$ for some $P \in \C^{\infty}(\R^3)$. However, if we take for $a_{13}$ a function only of $x_1$ and $x_2$ then then $\pl 3$ will commute with $\pl 2 +a_{13}\pl 3$. We can then take any $P_1 \in \C^{\infty}(\R^3)$ and set $P = \pl 3P_1$. Applying Theorem \ref{genODE} then gives a solution.

Applying Theorem \ref{varphi} to the Poisson structure on the dual of the simple  three dimensional Lie algebra, cf \S \ref{Poissonstructure}, also gives a large family of integrable infinitesimal deformations on $\C^3$. Theorem \ref{varphi} gives a partial order on integrable infinitesimals of $C^{\infty}(\R^3)$, with $F_1 \prec F_1' $ whenever $F_1' = \varphi F_1$ for some $\varphi$, and an equivalence relation when also $F_ 1\prec F_1' \prec F_1$.  It also implies, for dimension 3, that if $F_1$ is a basic infinitesimal of $C^{\infty}(\R^3)$ then $\varphi F_1$ is integrable for any $\varphi$ in $C^{\infty}(\R^3)$; such infinitesimals will be called quasibasic. 
For a basic infinitesimal $F_1 = f_1f_2\pl1\cpr\pl2+f_1f_3\pl1\cpr\pl3 + f_2f_3\pl2\cpr\pl3$ the stronger equations
\begin{equation}\label{separated}
\pl1(a_{12}/a_{13})  = \pl3(a_{13}/a_{23})  = \pl2(a_{23}/a_{12}) 
 = 0 \end{equation}
 hold. The analog is true in all dimensions.

\begin{theorem}\label{ind}
Let $F_1 =  \sum_{i<j}a_{ij}\pl i\cpr\pl j$  be an infinitesimal deformation of $\Cinf, $ where $d \ge 3$.  Suppose that the $a_{ij}$ are all invertible, and  that for all  
$\{i,j,k\}$ with $i<j<k$ we have $\pl i(a_{ij}/a_{ik}) = \pl j(a_{ij}/a_{jk})  =  \pl k(a_{ik}/a_{jk}) = 0$. Denote  \linebreak $a_{ij}(x_1,\dots,x_{i-1}, 0, x_{i+1},\dots, x_d)$  by $a_{ij}(\hat{x_i})$. Then $a_{ij}/a_{ij}(\hat{x_i})$ can be written in the form $1+x_iu_i$, where $u_i$ is not a function of $x_i$. 
Further, if the dimension is 3 then
$F_1$ is quasibasic; there are single variable functions $f_1,f_2,f_3$ of $x_1,x_2,x_3$, respectively, and a function $\varphi$ such that $$ a_{12} = \varphi f_1f_2, \quad a_{13} = \varphi f_1f_3\quad a_{23} = \varphi f_2f_3.$$
\end{theorem}
\noindent\textsc{Proof.}
For the first part, we can write $a_{ij}/a_{ij}(\hat{x_i}) = 1+x_iu_{ij}$ for certain functions $u_{ij}$. Suppose that $i<j<k$; the other cases are similar. By the hypotheses,
$[a_{ij}/a_{ij}(\hat{x_i})]/[a_{ik}/a_{ik}(\hat{x_i})]$ is not a function of $x_i$ but equals $(1+x_iu_{ij})/(1+x_iu_{ik})$,  so the latter is not a function of $x_i$. However, when $x_i = 0$ it is equal to 1, so in fact $(1+x_iu_{ij})/(1+x_iu_{ik})$ is identically equal to 1, implying that $u_{ij} = u_{ik}$. We may therefore write simply $u_i$ for $u_{ij}$.

Now
let $d=3$, so 
$F_1 = a_{12}\pl1\cpr\pl2 + a_{13}\pl1\cpr\pl3 +a_{23}\pl2\cpr\pl3. $
We may assume, without loss of generality, that $a_{12} = 1$. From the hypotheses it then follows that $a_{13}$ is not a function of $x_1$ and $a_{23}$ is not a function of $x_2$. It follows that $a_{13}(\hat{x_3})$ is a function only of $x_2$, which we write as $1/f_2$,  and that $a_{23}(\hat{x_3})$ is a function only of $x_1$, which we write as $1/f_1$. Then $a_{13} = 
(1+x_3u_3)/f_2$, where $u_3$ can be a function only of $x_2$ and $x_3$, and similarly $a_{23} = (1+x_3u_3)/f_2$, where $u_3$ now can be a function only of $x_1$ and $x_3$. Therefore, it is a function of $x_3$ alone. Denoting $(1+x_3u_3)$ by $f_3$, the original $F_1$ is therefore equivalent to an infinitesimal with $a_{12} = 1, a_{13} = f_3/f_2, a_{23} = f_3/f_1$. Multiplying by $f_1f_2$ gives the desired result.
$\Box$\medskip

We do not know if the analogous result holds in higher dimensions.

\section{The Big Bang}\label{Big Bang}
Theorem \ref{varphi} indicates an important property of the exceptional dimension 3, which is the smallest in which obstructions can appear. It allows the construction of infinitesimal deformations of $\C^{\infty}(\R^3)$ which are arbitrarily large somewhere inside a bounded  open set $U$ but which vanish outside $U$. For example, 
let $\Pi$ be a non-zero Poisson structure on $\C^{\infty}(\R^3)$, e.g., the Kirillov-Kostant-Souriau Poisson structure on the dual of the Lie algebra ${sl}_2$ of $2 \times 2$ matrices of trace zero, and $\varphi$ be a 
 smooth function which is arbitrarily large in a smaller open set and which vanishes outside $U$; take the infinitesimal to be $\varphi \Pi$.
 
The question now arises of whether there is an integral of the foregoing $\varphi \Pi$ to a full deformation of $\C^{\infty}(\R^3)$ which shares the property of vanishing outside $U$ and being arbitrarily large somewhere inside. The integral provided by Kontsevich's UDF has this property but is purely formal; to be physically meaningful, the integral must in addition be effective. If such exists, then by taking $U$ arbitrarily small one can have a deformation of $\C^{\infty}(\R^3)$ where outside $U$ the products of smooth functions are unchanged but inside this arbitrarily small $U$ the deformation is large. The problem is that the chosen $\varphi$ can not be analytic, hence neither is $\varphi \Pi$, and Kontsevich's UDF is not effective for such infinitesimal deformations. However, as an approximation, take $\varphi$ to be a highly peaked Gaussian distribution which becomes very small outside $U$. The effect of the infinitesimal deformation $\varphi \Pi$ may be said to be largely confined to $U$. It is analytic and one can apply Kontsevich's UDF.  

Now let  
 $\psi(t)$ be a smooth function of time. Then $\psi(t)\varphi(x_1,x_2,x_3)\Pi$ will again be integrable as an infinitesimal deformation of 
$ C^{\infty}(x_1,x_2,x_3,t)$ since $\Pi$ involved only the space variables. We may take $\psi(t)$ to be arbitrarily large at some point $t_0$ but vanishing outside an arbitrarily small neighborhood of $t_0$.
 Taking increasingly peaked Gaussian distributions and smaller intervals of time,
 in the limit, if that were meaningful, one would have an extreme disturbance localized at a point in space and instant of  time. 
 This might yield some insight into
 the nature of three-dimensional space and the Big Bang.


\begin{thebibliography}{50}

\bibitem{BFFLS}
F.~Bayen, M.~Flato, C.~Fr{\o}nsdal, A.~Lichnerowicz, and D.~Sternheimer.
\newblock {Deformation theory and quantization, I and II.}
\newblock{\em Ann. of Phys.}, 111:61--151, 1978.

\bibitem{CalaqueVandenBergh}
D.~Calaque and M.~Van~den~Bergh.
\newblock{Hochschild cohomology and Atiyah classes}.
\newblock{\em Adv. Math.} 224: 1839--1889, 2010 
\newblock{arXiv:0708.2725}.

\bibitem{CattFeldTom}
A.~Cattaneo,  D.~Felder, and L.~Tomassini.
\newblock{From local to global deformation quantization of Poisson manifolds.}
\newblock{\em Duke Math. J.} 115:2002, 329--352.

\bibitem{Coffee:Filtered}
Jane~Purcell~Coffee.
\newblock{Filtered and associated graded rings.}
\newblock{\em Bull. Amer. Math. Soc.} 78:1972, 584--587


\bibitem{CGG:noncomm}
V.~E.~Coll, M.~Gerstenhaber, A.~ Giaquinto, A.
\newblock{An explicit deformation formula with noncommuting  derivations}.
\newblock{in \textsc{Ring theory, Israel Math. Conf. Proc., 1}, Weizmann, Ramat Gan and Jerusalem, 1988/1989} 396--403.


\bibitem{DitoSternheimer:DefQuant}
G.~Dito and D.~Sernheimer.
\newblock{Deformation quantization: genesis, developments and metamorphoses}
\newblock{in \textsc{Deformation Quantization. 
Proceedings of the meeting between mathematicians and theoretical physicists, Strasbourg, 2001. G. Halbout, ed. IRMA Lectures in Math. Theoret. Phys., vol. 1}, Walter De Gruyter, Berlin, 2002} 9--54;
\newblock{also in arXiv:math/0201168v1 [math.QA]}.


\bibitem{Dolg:rational}
V.~Dolgushev.
\newblock{
A formality quasi-isomorphism for Hochschild cochains over rationals can be constructed recursively}.
\newblock{arXiv:1306.6733v4 [math.KT]  9 Feb 2017}.



\bibitem{Douady:Obstruction}
A.~Douady.
\raggedright \newblock{Obstruction primaire {\`a} la d{\'e}formation}.
\newblock{S{\'e}minaire Henri Cartan, tome 13, n$^o 1$, 1960-1961, 
exp. n$^o 4$, p. 1-19}.
\newblock{available at \linebreak
 \verb! http://numdam.org/item?id=SHC_1960-1961__13_1_A3_0!}.


\bibitem{FN}
A.~Froelicher and A.~Nijenhuis.
\newblock {A theorem on stability of complex structures}.
\newblock {\em Proc. Nat. Acad. Sci., U.S.A}, 43:239--241, 1957.

\bibitem{G:Discontinuous}
M.~Gerstenhaber
\newblock{On the algebraic structure of discontinuous groups}.
\newblock{\em Proc. Amr. Math. Soc.} 4:745--750, 1954.

\bibitem{G:Cohomology}
M.~Gerstenhaber.
\newblock{The cohomology structure of an associative ring}.
\newblock{\em Ann. of Math.} 78:267--288, 1963.

\bibitem{G:Uniform}
M.~Gerstenhaber.
\newblock{A uniform cohomology theory for algebras}.
\newblock{\em Proc. Nat. Acad. Sci. USA.} 51:626--629, 1964.

\bibitem{G:Def}
M.~Gerstenhaber.
\newblock {On the deformation of rings and algebras}.
\newblock {\em Ann. of Math.}, 79:59--103, 1964.

\bibitem{G:Def2}
M.~Gerstenhaber.
\newblock {On the deformation of rings and algebras: II}.
\newblock {\em Ann. of Math.}, 84:1 -- 19, 1966.



\bibitem{G:Def3}
M.~Gerstenhaber.
\newblock {On the deformation of rings and algebras: III }.
\newblock {\em Ann. of Math.}, 88:1--34, 1968.

\bibitem{G:Def4}
M.~Gerstenhaber.
\newblock {On the deformation of rings and algebras: IV }.
\newblock {\em Ann. of Math.}, 99:257--276, 1974.

\bibitem{G:PathAlgebras}
M.~Gerstenhaber.
\newblock{Path algebras, wave-particle duality, and quantization of phase space.} 
\newblock{\em Lett. Math. Phys.} 107:409--426, 2017.

\bibitem{GG:qWeyl}
M.~Gerstenhaber and A.~Giauinto.
\newblock{On the cohomology of the Weyl algebra, the quantum plane, and the $q$-Weyl algebra}.
\newblock{\em J. Pure and Applied Algebra}, 218:879--887, 2014.

Compatible deformations
Gerstenhaber, Murray; Giaquinto, Anthony; Green, Edward L; Huisgen-Zimmermann, Birge
ISSN: 0271-4132 , 1098-3627; ISBN: 0821809288 , 9780821878200; DOI: 10.1090/conm/229/03328
Trends in the representation theory of finite dimensional algebras : 1997 joint summer research conference on trends in the representation theory of finite dimensional algebras, July 20-24, 1997, Seattle, Washington / , 1998, Vol.229, p.159-168

\bibitem{GS:Hodge}
M.~Gerstenhaber and S.~D.~Schack.
\newblock{ A Hodge-type decomposition for commutative algebra cohomology}. 
\newblock{\em J. Pure Appl. Algebra}  48:229--247, 1987.

\bibitem{GS:AlgebraicCohomology}
M.~Gerstenhaber and S.~D.~Schack.
\newblock{Algebraic cohomology and deformation theory}.
\newblock{in \textsc{Deformation Theory of Algebras and Structures and
  Applications}, volume 247 of {\em NATO ASI Science Series}}.
\newblock Kluwer Academic Publishers, Dordrecht/Boston/London, 1988, 11 -- 264.

\bibitem{GiaqZhang:UDF} 
A.~Giaquinto and J.J.~Zhang.
\newblock{ Bialgebra actions, twists, and universal deformation formulas}. 
\newblock{\em J. Pure Appl. Algebra no. 2} 128:133--151, 1998

\bibitem{Groenewold}
H.~J. Groenewold.
\newblock {On the principles of elementary quantum mechanics}.
\newblock {\em Physica}, 12:405--460, 1946.

\bibitem{Heisenberg}
W.~Heisenberg.
\newblock{Ueber den anschaulichen Inhalt der quantentheoretischen Kinematik and Mechanik}.
\newblock{\em Zeitschrift für Physik}, 43: 172--198, 1927.
\newblock{English translation in} 
Wheeler, J.~A. and W.~H.~Zurek (eds), 1983, {\sc Quantum Theory and Measurement}, Princeton, NJ: Princeton University Press,   62--84, 1983.

\bibitem{HKR}
G.~Hochschild, B.~Kostant, A.~Rosenberg.
\newblock {Differential forms on regular affine algebras}.
\newblock {\em Trans. Amer. Math. Soc.}, 102:383--408, 1962.

\bibitem{Kodaira:Complex}
K.~Kodaira.
 \newblock{\textsc {Complex Manifolds and Deformation of Complex Structures}}.
\newblock{Classics in Mathematics, Springer, Berlin Heidelberg New York, 2005}.
\newblock{Originally published as vol. 283 in the series \em{Grundlehren der mathematischen Wissenschaften}}.

\bibitem{Kontsevich:Poisson}
M.~Kontsevich.
\newblock{Deformation quantization of Poisson manifolds}.
\newblock{\em Lett. Math. Phys.} 66:157--216, 2003.
\newblock{arXiv:9709040 [math.QA]}.

\bibitem{Mac Lane:Homology}
S.~Mac~Lane.
\newblock{\textsc{Homology}, Third Corrected Printing, Springer, New York Heidelberg Berlin, 1975}.

\bibitem{Matt:Exp}
S.~Mattarei
\newblock{Exponential functions in prime characteristic}.
\newblock{arXiv:math/0511168v1 [math.GM]}.
\newblock{7 Nov 2005; accessed 20 Oct 2022}.


\bibitem{Moyal}
J.~E.~Moyal.
\newblock {Quantum mechanics as a statistical theory}.
\newblock {{Proc. Cambridge Phil. Soc.}} 45:99 --124, 1949.


\bibitem{Roger}
C.~Roger.
\newblock {Gerstenhaber and Batalin-Vilkovisky algebras; algebraic, geometric, and physical aspects}.
\newblock{\em Archivum Mathematicum (Brno)}, 45:301--324, 2009.

\bibitem{Schiller}
J.~J.~Schiller.
\newblock{Moduli for Special Riemann Surfaces of Genus 2}.
\newblock{\em Trans. Amer. Math. Soc.Trans. Amer. Math. Soc.}, 144:95--113,1969.


\bibitem{SchlessStash}
M.~Schlessinger and J. Stasheff.
\newblock{Deformation theory and rational homotopy type}.
\newblock{arXiv:1211.1647v1 [math.QA]} 7 Nov 2012.

\bibitem{Srid:thesis}
R.~Sridharan.
\newblock{Filtered algebras and representations of Lie algebras}.
\newblock{\em Trans. Amer. Math. Soc.}, 100:530–550, 1961. 

\bibitem{Stasheff:Hspaces}
J.~D.~Stasheff.
\newblock{Homotopy associativity of H-spaces. I, II}. 
\newblock{ \em Trans. Amer. Math. Soc.} 108, 275-292; ibid.  293–312, 1963.

\bibitem{Sternheimer:20YearsAfter}
D.~Sternheimer.
\newblock{Deformation quantization: Twenty years after}, in \textsc{Particles, Fields, and Gravitation}, {\L}odz, 1998 (J. Rembieli{\'n}ski, ed.)
\newblock AIP Press, New York, 1998, 107--145.
\newblock {arXiv:9809056 [math.QA] 10 Sep 1998}.

\bibitem{Tamarkin:formality}
D.~E.~Tamarkin. 
\newblock{Another proof of Kontsevich formality conjecture}.
\newblock{arXiv:9803025 [math.QA] 24 Sep 1998}.

\bibitem{Teich:Extremale}
O.~Teichm\"uller.
\newblock{Extremale quasikonforme Abbildungen und quadratische Differentiale}.
\newblock {\em Abhandlungen Preussische Akad. Wiss. Math.-Nat. Kl.} 1939 (no. 22):197 pp, 1940.
\newblock see also {\em Oswald Teichm\"uller: Gesammelete Abhandlungen}.
\newblock {Lars. V. Ahlfors and Frederick W. Gehring, ed}.
\newblock {Springer-Verlag, Berlin-New York, 1982}.


\end{thebibliography}
\end{document}